\documentclass[12pt]{article}
\usepackage{graphicx}
\usepackage{amsmath,amsthm,amssymb,enumerate}%, esint}
\usepackage{euscript,mathrsfs}
\usepackage{color}
\usepackage[left=2cm,right=2cm,top=3.5cm,bottom=3.5cm]{geometry}
\usepackage{color}
\catcode`\@=11 \@addtoreset{equation}{section}

\catcode`\@=12

\allowdisplaybreaks

\newtheorem{Theorem}{Theorem}[section]
\newtheorem{Proposition}[Theorem]{Proposition}
\newtheorem{Lemma}[Theorem]{Lemma}
\newtheorem{Corollary}[Theorem]{Corollary}

\theoremstyle{definition}
\newtheorem{Definition}[Theorem]{Definition}

\newtheorem{Remark}[Theorem]{Remark}

\newcommand{\bTheorem}[1]{
\begin{Theorem} \label{T#1} }
\newcommand{\eT}{\end{Theorem}}

\newcommand{\bProposition}[1]{
\begin{Proposition} \label{P#1}}
\newcommand{\eP}{\end{Proposition}}

\newcommand{\bLemma}[1]{
\begin{Lemma} \label{L#1} }
\newcommand{\eL}{\end{Lemma}}

\newcommand{\bCorollary}[1]{
\begin{Corollary} \label{C#1} }
\newcommand{\eC}{\end{Corollary}}

\newcommand{\bRemark}[1]{
\begin{Remark} \label{R#1} }
\newcommand{\eR}{\end{Remark}}

\newcommand{\bDefinition}[1]{
\begin{Definition} \label{D#1} }
\newcommand{\eD}{\end{Definition}}

\newcommand{\bfphi}{\boldsymbol{\varphi}}

\newcommand{\bFormula}[1]{
\begin{equation} \label{#1}}
\newcommand{\eF}{\end{equation}}

\newcommand{\Ov}[1]{\overline{#1}}

\newcommand{\DC}{C^\infty_c}

\newcommand{\vr}{\varrho}

\newcommand{\vu}{\vc{u}}

\newcommand{\vc}[1]{{\bf #1}}

\newcommand{\Div}{{\rm div}_x}
\newcommand{\Grad}{\nabla_x}

\newcommand{\dx}{\,{\rm d} {x}}

\newcommand{\dt}{\,{\rm d} t }

\newcommand{\intO}[1]{\int_{\Omega} #1 \ \dx}

\definecolor{Cgrey}{rgb}{0.85,0.85,0.85}
\definecolor{Cblue}{rgb}{0.50,0.85,0.85}
\definecolor{Cred}{rgb}{1,0,0}
\definecolor{fancy}{rgb}{0.10,0.85,0.10}

\newcommand\Cbox[2]{%
    \newbox\contentbox%
    \newbox\bkgdbox%
    \setbox\contentbox\hbox to \hsize{%
        \vtop{
            \kern\columnsep
            \hbox to \hsize{%
                \kern\columnsep%
                \advance\hsize by -2\columnsep%
                \setlength{\textwidth}{\hsize}%
                \vbox{
                    \parskip=\baselineskip
                    \parindent=0bp
                    #2
                }%
                \kern\columnsep%
            }%
            \kern\columnsep%
        }%
    }%
    \setbox\bkgdbox\vbox{
        \color{#1}
        \hrule width  \wd\contentbox %
               height \ht\contentbox %
               depth  \dp\contentbox
        \color{black}
    }%
    \wd\bkgdbox=0bp%
    \vbox{\hbox to \hsize{\box\bkgdbox\box\contentbox}}%
    \vskip\baselineskip%
}

\date{}

\begin{document}

%%%%%%%%%%%%%%%%%%%%%%%%%%%%%%%%

\title{On weak--strong uniqueness for the compressible Navier--Stokes system 
with non--monotone pressure law}

\author{Eduard Feireisl\thanks{The research of E.F.~leading to these results has received funding from the
Czech Sciences Foundation (GA\v CR), Grant Agreement
18--05974S. The Institute of Mathematics of the Academy of Sciences of
the Czech Republic is supported by RVO:67985840.}}

\date{\today}

\maketitle

\bigskip

\centerline{Institute of Mathematics of the Academy of Sciences of the Czech Republic}

\centerline{\v Zitn\' a 25, CZ-115 67 Praha 1, Czech Republic}

\centerline{and}

\centerline{Institut fuer Mathematik, Technische Universitaet Berlin}

\centerline{Strasse des 17. Juni 136, 
D - 10623 Berlin, Germany}

\bigskip

\begin{abstract}

We show the weak--strong uniqueness property for the compressible Navier--Stokes system with general non--monotone pressure law. 
A weak solution coincides with the strong solution emanating from the same initial data as long as the latter solution exists.

\end{abstract}


%\tableofcontents

\section{Introduction}
\label{M}

Let $\Omega \subset R^N$, $N=1,2,3$ be a bounded Lipschitz domain. 
The Navier--Stokes system desribing 
the time evolution of the density $\vr = \vr(t,x)$ and the velocity 
$\vu = \vu(t,x)$ of a compressible barotropic viscous fluid reads:
\begin{equation} \label{M1}
\partial_t \vr + \Div (\vr \vu)  = 0,
\end{equation}
\begin{equation} \label{M2}
\begin{split}
\partial_t (\vr \vu) + \Div (\vr \vu \otimes \vu) &+ \Grad p(\vr)  = \Div \mathbb{S}(\Grad \vu), 
\end{split}
\end{equation}
where the viscous stress is given by Newton's rheological law 
\begin{equation} \label{M4}
\mathbb{S}(\Grad \vu) = \mu \left( \Grad \vu + \Grad^t \vu - \frac{2}{N} \Div \vu \mathbb{I} \right) + 
\lambda \Div \vu \mathbb{I}, \ \mu > 0, \ \lambda \geq 0.
\end{equation}
We consider the no--slip boundary condition
\begin{equation} \label{M3}
\vu|_{\partial \Omega} = 0,
\end{equation}
and the barotropic pressure law
\begin{equation} \label{M5}
p(\vr) = a \vr^\gamma + q(\vr), \ q \in \DC(0, \infty), \ a > 0, \ \gamma > 1.
\end{equation}

If $q \equiv 0$, the relation (\ref{M5}) reduces to the standard \emph{isentropic} equation of state, for which the problem 
(\ref{M1}--\ref{M5}) admits global in time weak solutions for any finite energy initial data, see Antontsev et al. 
\cite{AKM} for $N = 1$, Lions \cite{LI4} for $N = 2$, $\gamma \geq \frac{3}{2}$, $N=3$, $\gamma \geq \frac{9}{5}$, and 
\cite{EF70} for $N = 2$, $\gamma > 1$, $N = 3$, $\gamma > \frac{3}{2}$. 

If $q \ne 0$, the pressure $p$ need not be a monotone function of the density. The weak solutions, however, still exist globally in time, at least for $\gamma > \frac{3}{2}$ and $N = 3$, see \cite{EF61}. The result has been extended recently to more general 
(not necessarily compactly supported)
perturbations $q$ and $\gamma \geq 2$ by Bresch and Jabin \cite{BreJab}.

If the initial data are smooth enough, the same problem admits local in time strong solutions that are global if $N=1$ or $N=2,3$ and the data are sufficiently small, see 
\cite{AKM}, Matsumura and Nishida \cite{MANI}, among others. A natural question arises whether strong solutions are uniquely determined in the class of weak solutions - a weak solution and the strong solution starting from the same initial data coincide on the life 
span of the latter. The first result of this type was shown by Germain \cite{Ger} in the class of weak solutions enjoying certain additional regularity properties. Finally, the weak--strong uniqueness property was established in \cite{FeJiNo}, \cite{FENOSU} 
in the class of dissipative weak solutions, the existence of which is guaranteed by the above mentioned existence theory. 

The weak--strong uniqueness property is strongly related to the convexity of the energy functional 
\[
[\vr, \vc{m}] \mapsto \frac{1}{2} \frac{|\vc{m}|^2}{\vr} + H(\vr),\ 
H(\vr) = \vr \int_1^\vr \frac{ p(z) }{z^2} \ {\rm d}z. 
\]
In particular, the pressure $p(\vr)$ must be (strictly) increasing function of $\vr$ as $H''(\vr) = p'(\vr)/ \vr$, which excludes 
the general pressure law (\ref{M5}) with $q \not\equiv 0$. The goal of this short note is to show that the technique of \cite{FeJiNo}
can be accommodated to handle a general non--monotone pressure law (\ref{M5}).     

\section{Dissipative weak solutions, main result}
 
Suppose that $\gamma > 1$, $N=1,2,3$. We say that 
\[
\vr \in L^\infty (0,T; L^\gamma(\Omega)), \ \vr \geq 0,\ 
\vu \in L^2(0,T; W^{1,2}_0 (\Omega; R^N)),\ 
\vc{m} \equiv \vr \vu \in L^\infty(0,T; L^{\frac{2 \gamma}{\gamma + 1}} (\Omega; R^N)),
\]
is a \emph{dissipative weak solution} to problem (\ref{M1}--\ref{M5}) if:
\begin{itemize}
\item the integral identity 
\begin{equation} \label{D1}
\left[ \intO{ \vr \varphi } \right]_{t = 0}^{t = \tau} = 
\int_0^\tau \intO{ \left[ \vr \partial_t \varphi + \vr \vu \cdot \Grad \varphi \right] } \dt
\end{equation}
holds for any $\tau \in [0,T]$ and any $\varphi \in C^1(\Ov{\Omega} \times [0,T])$; 
\item
the integral identity 
\begin{equation} \label{D2}
\begin{split}
&\left[ \intO{ \vr \vu \cdot \bfphi } \right]_{t = 0}^{t = \tau} 
\\
&= 
\int_0^\tau \intO{ \left[ \vr \vu \cdot \partial_t \bfphi + \vr \vu \otimes \vu : \Grad \bfphi + 
p(\vr) \Div \bfphi \right] } \dt - \int_0^\tau \intO{ \mathbb{S}(\Grad \vu) : \Grad \bfphi } \dt 
\end{split}
\end{equation}
holds
for any $\tau \in [0,T]$ and any $\bfphi \in C^1(\Ov{\Omega} \times [0,T]; R^N)$, $\bfphi|_{\partial \Omega} = 0$;
\item
the renormalized equation of continuity holds, meaning, the integral identity
\begin{equation} \label{D3}
\left[ \intO{ b(\vr) \varphi } \right]_{t = 0}^{t = \tau} = 
\int_0^\tau \intO{ \Big[ b(\vr) \partial_t \varphi + b(\vr) \vu \cdot \Grad \varphi 
+(b(\vr) - b'(\vr) \vr) \Div \vu \varphi \Big] } \dt
\end{equation}
holds
for any $\tau \in [0,T]$, $\varphi \in C^1(\Ov{\Omega} \times [0,T])$, and $b \in C^1[0, \infty)$, $b' \in C_c[0, \infty)$; 
\item 
the energy inequality 
\begin{equation} \label{D4}
\begin{split}
\left[ \intO{ \left( \vr |\vu|^2 + P(\vr) \right) } \right]_{t = 0}^{t=\tau} &+ 
\int_0^\tau \intO{ \mathbb{S}(\Grad \vu) : \Grad \vu } \dt \leq 0, \\
P(\vr) &= H(\vr) + Q(\vr), \ H(\vr) = \frac{a}{\gamma - 1} \vr^\gamma,\ 
Q(\vr) = \vr \int_1^{\vr} \frac{q(z)}{z^2} \ {\rm dz} 
\end{split}
\end{equation} 
holds for a.a. $\tau \in [0,T]$.

\end{itemize}

As $\vr$ satisfies (\ref{D1}), (\ref{D3}), we get 
\[
\left[ \intO{ Q(\vr) } \right]_{t = 0}^{t = \tau} = 
- \int_0^\tau \intO{ q(\vr) \Div \vu } \dt;
\] 
whence it follows from (\ref{D4}) that 
\begin{equation} \label{D5}
\left[ \intO{ \left( \vr |\vu|^2 + H(\vr) \right) } \right]_{t = 0}^{t=\tau} + 
\int_0^\tau \intO{ \mathbb{S}(\Grad \vu) : \Grad \vu } \dt \leq \int_0^\tau \intO{ q(\vr) \Div \vu } \dt. 
\end{equation} 
Relation (\ref{D5}) holds for \emph{any} $t \in [0,T]$ due to the weak lower semi--continuity of the functional
\[
[\vr, \vc{m} = \vr \vu] \mapsto \intO{ \left( \frac{|\vc{m}|^2}{\vr} + H(\vr) \right) }.
\]

Our goal is to show the following result. 

\begin{Theorem} \label{TD1}
Let $\Omega \subset R^N$ be a bounded Lipschitz domain. Let the pressure $p = p(\vr)$ be given by (\ref{M5}). 
Suppose that $[\vr, \vu]$ is a dissipative weak solution and $[r, \vc{U}]$ a classical solution of the problem 
(\ref{M1}--\ref{M5})
on the time interval $[0,T]$ such that 
\[
\vr(0, \cdot) = r(0, \cdot) > 0, \ \vr \vu(0, \cdot) = r(0, \cdot) \vc{U}(0, \cdot).
\] 

Then 
\[
\vr = r, \ \vu = \vc{U} \ \mbox{in}\ (0,T) \times \Omega.
\]

\end{Theorem}

The rest of the paper is devoted to the proof of Theorem \ref{TD1}. 

\section{Relative energy}

Following \cite{FeJiNo} (cf. the standard reference material by Dafermos \cite{Daf4})
we introduce the \emph{relative energy functional}:
\[
\begin{split}
\mathcal{E} \left( \vr, \vu \ \Big| \ r , \vc{U} \right) 
= \intO{ \left[ \frac{1}{2} \vr |\vc{u} - \vc{U}|^2 + H(\vr) - H'(r)(\vr - r) - H(r) \right] } 
= \sum_{j = 1}^4 I_j,
\end{split}
\]
where 
\[
I_1 = \intO{ \left( \frac{1}{2} \vr |\vu|^2 + H(\vr) \right) }, 
\]
\[
I_2 = - \intO{ \vr \vu \cdot \vc{U} } , \ I_3 = \intO{ \vr \left( \frac{1}{2} |\vc{U}|^2 - H'(r) \right) }
\]
\[
I_4 = \intO{ \left( H'(r) r - H(r) \right) } = \intO{ a r^\gamma }.
\]
Note that $\mathcal{E} \left( \vr, \vu \ \Big| \ r , \vc{U} \right)$ is well defined as soon as $[\vr, \vu]$ is a dissipative weak solution and $R$ and $\vc{U}$ are arbitrary continuous differentiable functions satisfying the natural compatibility conditions 
\begin{equation} \label{comp}
r \in C^1([0,T] \times \Ov{\Omega}),\ r > 0, \ \vc{U} \in C^1([0,T] \times \Ov{\Omega}; R^N),\ 
\vc{U}|_{\partial \Omega}  = 0.
\end{equation}
 
Using the weak formulation (\ref{D1}--\ref{D5}) we deduce easily 
\begin{equation} \label{R1}
\left[ \mathcal{E} \left( \vr, \vu \ \Big| \ r , \vc{U} \right) \right]_{t = 0}^{t = \tau} = 
\sum_{j=1}^4 [I_j ]_{t = 0}^{t = \tau},
\end{equation}
where 
\begin{equation} \label{R2}
[I_1]_{t = 0}^{t = \tau} + \int_0^\tau \intO{ \mathbb{S}(\Grad \vu) : \Grad \vu } \dt \leq \int_0^\tau \intO{ q(\vr) \Div \vu } \dt,
\end{equation}
\begin{equation} \label{R3}
\begin{split}
&[I_2]_{t = 0}^{t = \tau}\\ 
&= 
-\int_0^\tau \intO{ \left[ \vr \vu \cdot \partial_t \vc{U} + \vr \vu \otimes \vu : \Grad \vc{U} + 
p(\vr) \Div \vc{U} \right] } \dt + \int_0^\tau \intO{ \mathbb{S}(\Grad \vu) : \Grad \vc{U} } \dt\\
&= -\int_0^\tau \intO{ \left[ \vr \vu \cdot \partial_t \vc{U} + \vr \vu \otimes \vu : \Grad \vc{U} + 
a \vr^\gamma \Div \vc{U}  + q(\vr) \Div \vc{U} \right] } \dt\\
&+ \int_0^\tau \intO{ \mathbb{S}(\Grad \vu) : \Grad \vc{U} } \dt,
\end{split}
\end{equation}
and
\begin{equation} \label{R4}
\begin{split}
[I_3]_{t = 0}^{t = \tau} &= \int_0^\tau \intO{ \left[
\vr \vc{U} \cdot \partial_t \vc{U} + \vr \vu \cdot \vc{U} \cdot \Grad \vc{U} \right] } \dt \\ 
&- \int_0^\tau \intO{ \left[
\vr \vc{U} \cdot \partial_t H'(r) + \vr \vu \cdot \Grad H'(r) \right] } \dt, 
\end{split}
\end{equation}
cf. \cite{FeJiNo}. 

Summing up (\ref{R2}--\ref{R4}) we obtain the relative energy inequality
\begin{equation} \label{R5}
\begin{split}
&\left[ \mathcal{E} \left( \vr, \vu \ \Big| \ r , \vc{U} \right) \right]_{t = 0}^{t = \tau} + 
\int_0^\tau \intO{ \mathbb{S}(\Grad \vu): (\Grad \vu - \Grad \vc{U} ) } \dt \\
&\leq \int_0^\tau \intO{ \left[ \vr (\vc{U} - \vu)  \cdot \partial_t \vc{U} + \vr \vu \cdot (\vc{U} - \vu) \cdot \Grad \vc{U} - 
a \vr^\gamma \Div \vc{U}\right] } \dt\\
&+ \int_0^\tau \intO{ q(\vr) \left( \Div \vu - \Div \vc{U} \right) } \dt \\
&- \int_0^\tau \intO{ \left[
\vr \vc{U} \cdot \partial_t H'(r) + \vr \vu \cdot \Grad H'(r) \right] } \dt + 
\int_0^\tau \intO{ a \partial_t r^\gamma } \dt 
\end{split}
\end{equation}
for any $\tau \in [0,T]$ and any $r$ and $\vc{U}$ satisfying (\ref{comp}).  

\section{Weak strong uniqueness}

We show Theorem \ref{TD1} by considering 
the strong solution $[r, \vc{U}]$ as test functions in the relative energy inequality (\ref{R5}). 

\medskip

\noindent
$\bullet$ {\bf Step 1} 

We write 
\[
\intO{ \vr \vu \cdot (\vc{U} - \vu) \cdot \Grad \vc{U} } 
= \intO{ \vr (\vu - \vc{U}) \cdot (\vc{U} - \vu) \cdot \Grad \vc{U} }
+ \intO{ \vr \cdot (\vc{U} - \vu) \cdot \vc{U} \cdot \Grad \vc{U} }
\]
where
\[
\intO{ \vr (\vu - \vc{U}) \cdot (\vc{U} - \vu) \cdot \Grad \vc{U} } 
\leq c_1 \mathcal{E} \left( \vr, \vu \ \Big| \ r , \vc{U} \right). 
\]
As 
\[
\partial_t \vc{U} + \vc{U} \cdot \Grad \vc{U} = - \frac{1}{r} \Grad p(r) + \frac{1}{r} \Div \mathbb{S}(\Grad \vc{U}) 
\]
we deduce from (\ref{R5})
\begin{equation} \label{W1}
\begin{split}
&\left[ \mathcal{E} \left( \vr, \vu \ \Big| \ r , \vc{U} \right) \right]_{t = 0}^{t = \tau} + 
\int_0^\tau \intO{ \mathbb{S}(\Grad \vu): (\Grad \vu - \Grad \vc{U} ) } \dt \\
&\leq \int_0^\tau \intO{ \left[ \frac{\vr}{r} (\vc{U} - \vu) \Div \mathbb{S}(\Grad \vc{U})  - \frac{\vr}{r} (\vc{U} - \vu) \cdot 
\Grad p(r)  - 
a \vr^\gamma \Div \vc{U}\right] } \dt\\
&+ \int_0^\tau \intO{ q(\vr) \left( \Div \vu - \Div \vc{U} \right) } \dt \\
&- \int_0^\tau \intO{ \left[
\vr \vc{U} \cdot \partial_t H'(r) + \vr \vu \cdot \Grad H'(r) \right] } \dt + 
\int_0^\tau \intO{ a \partial_t r^\gamma } \dt \\
&+ c_1 \int_0^\tau \intO{ \mathcal{E} \left( \vr, \vu \ \Big| \ r , \vc{U} \right) } \dt 
\end{split}
\end{equation}

\medskip

\noindent $\bullet$ {\bf Step 2}

Using the relation $p(r) = a r^\gamma + q(r)$ we may regroup terms in (\ref{W1}) obtaining
\[
\begin{split}
&\left[ \mathcal{E} \left( \vr, \vu \ \Big| \ r , \vc{U} \right) \right]_{t = 0}^{t = \tau} + 
\int_0^\tau \intO{ \left( \mathbb{S}(\Grad \vu) - \mathbb{S}(\Grad \vc{U}) \right) : (\Grad \vu - \Grad \vc{U} ) } \dt \\
&\leq \int_0^\tau \intO{ \left[ \left( \frac{\vr}{r} - 1 \right) (\vc{U} - \vu) \Div \mathbb{S}(\Grad \vc{U})  - a \frac{\vr}{r} (\vc{U} - \vu) \cdot 
\Grad r^\gamma   - 
a \vr^\gamma \Div \vc{U}\right] } \dt\\
& + \int_0^\tau \intO{ \left( \frac{\vr}{r} - 1 \right) (\vu - \vc{U}) \cdot 
\Grad q(r) } \dt  
\\
&+ \int_0^\tau \intO{ \Big( q(\vr) - q(r) \Big)  \left( \Div \vu - \Div \vc{U} \right)  } \dt \\
&- \int_0^\tau \intO{ \left[
\vr \vc{U} \cdot \partial_t H'(r) + \vr \vu \cdot \Grad H'(r) \right] } \dt + 
\int_0^\tau \intO{ a \partial_t r^\gamma } \dt \\
&+ c_1 \int_0^\tau \intO{ \mathcal{E} \left( \vr, \vu \ \Big| \ r , \vc{U} \right) } \dt 
\end{split}
\]
As both $\vu$ and $\vc{U}$ satisfy the no--slip boundary conditions, we have 
\[
\left\| \Grad \vu - \Grad \vc{U} \right\|^2_{L^2(\Omega; R^N)} \leq c_2
\intO{ \left( \mathbb{S}(\Grad \vu) - \mathbb{S}(\Grad \vc{U}) \right) : (\Grad \vu - \Grad \vc{U} ) }, 
\]
and, consequently, 
\[
\begin{split}
&\intO{ \Big( q(\vr) - q(r) \Big)  \left( \Div \vu - \Div \vc{U} \right)  }\\  
&\leq c_3 \intO{ \Big( q(\vr) - q(r) \Big)^2 } + 
\frac{1}{2} \intO{ \left( \mathbb{S}(\Grad \vu) - \mathbb{S}(\Grad \vc{U}) \right) : (\Grad \vu - \Grad \vc{U} ) }.
\end{split}
\]

Thus we may infer that
\begin{equation} \label{W2}
\begin{split}
&\left[ \mathcal{E} \left( \vr, \vu \ \Big| \ r , \vc{U} \right) \right]_{t = 0}^{t = \tau} + 
\frac{1}{2} \int_0^\tau \intO{ \left( \mathbb{S}(\Grad \vu) - \mathbb{S}(\Grad \vc{U}) \right) : (\Grad \vu - \Grad \vc{U} ) } \dt \\
&\leq \int_0^\tau \intO{  \left( \frac{\vr}{r} - 1 \right) (\vc{U} - \vu) \cdot \Big( \Div \mathbb{S}(\Grad \vc{U}) 
- \Grad q(r) \Big)  } \dt\\
&+ c_4 \int_0^\tau \intO{ \Big( q(\vr) - q(r) \Big)^2} \dt \\
&- \int_0^\tau \intO{ \left[ a \frac{\vr}{r} (\vc{U} - \vu) \cdot 
\Grad r^\gamma  + a \vr^\gamma \Div \vc{U} \right] } \dt 
\\
&- \int_0^\tau \intO{ \left[
\vr \vc{U} \cdot \partial_t H'(r) + \vr \vu \cdot \Grad H'(r) \right] } \dt + 
\int_0^\tau \intO{ a \partial_t r^\gamma } \dt\\ &+ c_4 \int_0^\tau \intO{ \mathcal{E} \left( \vr, \vu \ \Big| \ r , \vc{U} \right) } \dt 
\end{split}
\end{equation}

\medskip

\noindent $\bullet$ {\bf Step 3}

Seeing that 
\[
H''(r) = a (\gamma - 1) r^{\gamma - 2}
\]
we obtain, after a simple manipulation for which we refer to \cite{FeJiNo},
\[
\begin{split}
&- \int_0^\tau \intO{ \left[ a \frac{\vr}{r} (\vc{U} - \vu) \cdot 
\Grad r^\gamma  + a \vr^\gamma \Div \vc{U} \right] } \dt 
\\
&- \int_0^\tau \intO{ \left[
\vr \vc{U} \cdot \partial_t H'(r) + \vr \vu \cdot \Grad H'(r) \right] } \dt + 
\int_0^\tau \intO{ a \partial_t r^\gamma } \dt 
\\
&= - \int_0^\tau \intO{ \Div \vc{U} \left( h(\vr) - h'(r)(\vr - r) - h(r) \right) } \dt 
\end{split}
\] 
where we have denoted $h(\vr) = a \vr^\gamma$. 

Consequently, (\ref{W2}) reduces to 
\begin{equation} \label{W3}
\begin{split}
&\left[ \mathcal{E} \left( \vr, \vu \ \Big| \ r , \vc{U} \right) \right]_{t = 0}^{t = \tau} + 
\frac{1}{2} \int_0^\tau \intO{ \left( \mathbb{S}(\Grad \vu) - \mathbb{S}(\Grad \vc{U}) \right) : (\Grad \vu - \Grad \vc{U} ) } \dt \\
&\leq \int_0^\tau \intO{  \left( \frac{\vr}{r} - 1 \right) (\vc{U} - \vu) \cdot \Big( \Div \mathbb{S}(\Grad \vc{U}) 
- \Grad q(r) \Big)  } \dt\\
&+ c_4 \int_0^\tau \intO{ \Big( q(\vr) - q(r) \Big)^2} \dt+ c_5 \int_0^\tau \intO{ \mathcal{E} \left( \vr, \vu \ \Big| \ r , \vc{U} \right) } \dt. 
\end{split}
\end{equation}

\medskip

\noindent
$\bullet$ {\bf Step 4}

Finally, we introduce a cut--off function $\Psi \in \DC(0, \infty)$, 
\[
0 \leq \Psi \leq 1, \ \Psi \equiv 1 \ \mbox{in} \ [\delta, \frac{1}{\delta}], 
\]
where $\delta$ is chosen so small that 
\[
r(t,x) \in [2 \delta, \frac{1}{2 \delta}]\ \mbox{for all}\ (t,x), \ {\rm supp}[q] \subset [2 \delta, \frac{1}{2 \delta}].
\]
Moreover, for $h \in L^1((0,T) \times \Omega)$, we set 
\[
h = h_{{\rm ess}} + h_{\rm res},\  h_{{\rm ess}} = \Psi (\vr) h,\ 
h_{\rm res} = (1 - \Psi (\vr)) h.
\]
It is easy to check that 
\[
\begin{split}
\mathcal{E} \left( \vr, \vu \ \Big| \ r , \vc{U} \right) \geq
c_6 \intO{ \left( [ \vu - \vc{U} ]^2_{\rm ess} + [ \vr - r ]^2_{\rm ess} + 1_{\rm res} + \vr^\gamma_{\rm res} \right)} .   
\end{split} 
\]

Consequently, we get 
\[
\begin{split}
\intO{ \Big( q(\vr) - q(r) \Big)^2} &\leq \intO{ \Big[ q(\vr) - q(r) \Big]^2_{\rm ess} } 
+ \intO{ \Big[ q(\vr) - q(r) \Big]^2_{\rm ess} } \\&\leq
c_7 \left[ \intO{ \Big[ \vr - r \Big]^2_{\rm ess} } + \intO{ q(r)^2_{\rm res} } \right] \leq 
c_9 \mathcal{E} \left( \vr, \vu \ \Big| \ r , \vc{U} \right)   
\end{split}
\]

Similarly, 
\[
\begin{split}
&\intO{  \left( \frac{\vr}{r} - 1 \right) (\vc{U} - \vu) \cdot \Big( \Div \mathbb{S}(\Grad \vc{U}) 
- \Grad q(r) \Big)  }\\
&\leq c_{10} \intO{ |\vr - r | \ |\vc{U} - \vu | } 
\leq c_{10} \left[ \intO{ | [\vr - r]_{\rm ess} | \ |\vc{U} - \vu | } 
+ \intO{ | [\vr - r]_{\rm res} | \ |\vc{U} - \vu | } \right] \\
&c_{11}(\delta) \left[ \mathcal{E} \left( \vr, \vu \ \Big| \ r , \vc{U} \right) + 
\intO{ \vr |\vu - \vc{U} |^2 } + \intO{ [1 + \vr]_{\rm res} } + \delta \| \vu - \vc{U} \|^2_{L^2(\Omega; R^N)} \right]     
\end{split}
\]
for any $\delta > 0$, where, by means of 
the Poincar\`e inequality, 
\[
\| \vu - \vc{U} \|^2_{L^2(\Omega; R^N)} \leq c_{12} \intO{ \left( \mathbb{S}(\Grad \vu) - \mathbb{S}(\Grad \vc{U}) \right) : (\Grad \vu - \Grad \vc{U} ) }. 
\]

Thus, 
going back to (\ref{W3}), we conclude 
\[
\left[ \mathcal{E} \left( \vr, \vu \ \Big| \ r , \vc{U} \right) \right]_{t = 0}^{t = \tau} 
\leq c_{13} \int_0^\tau \intO{ \mathcal{E} \left( \vr, \vu \ \Big| \ r , \vc{U} \right) } \dt. 
\]
Aplying Gronwall lemma we complete the proof of Theorem \ref{TD1}. 

\section{Concluding remarks}

The hypotheses concerning the pressure law can be relaxed, in particular, we may handle the pressure satisfying the hypotheses of 
\cite{EF61}. The result can be extended to the class of measure--valued solutions in the spirit of \cite{FGSWW1}. The method, however, cannot be extended to the Euler (inviscid) system as the presence of the viscous damping plays a crucial role in the proof.

\def\cprime{$'$} \def\ocirc#1{\ifmmode\setbox0=\hbox{$#1$}\dimen0=\ht0
  \advance\dimen0 by1pt\rlap{\hbox to\wd0{\hss\raise\dimen0
  \hbox{\hskip.2em$\scriptscriptstyle\circ$}\hss}}#1\else {\accent"17 #1}\fi}

%\bibliographystyle{plain}
%\bibliography{citace}

\end{document}